\documentclass[a4paper,twoside,11pt]{article}
\usepackage{amssymb}
\usepackage{amsmath}
\usepackage{amsthm}
\usepackage{latexsym}
\usepackage{amsfonts}
\usepackage{graphicx}
\usepackage{graphics}

\newcommand{\dis}{\displaystyle}
\theoremstyle{plain}
\newtheorem{thm}{Theorem}[section]   

\newtheorem{lem}[thm]{Lemma}

\newtheorem{Def}[thm]{Definition}

\theoremstyle{definition}
\newtheorem{rem}[thm]{Remark}

\newtheorem*{Proof}{Proof}

\newcommand{\bbb}[1]{\mbox{\boldmath$#1$}}

\newcommand{\el}{\ell}

\newcommand{\ra}{\;\rightarrow\;}

\newcommand{\ga}{\gamma }

\newcommand{\Sig} {{\varSigma}}
\newcommand{\de}{\delta }
\newcommand{\OO} {{\varOmega}}

\newcommand{\e}{\varepsilon }
\newcommand{\f}{\varphi}
\newcommand{\Fi}{\varPhi}

\newcommand{\zi}{\zeta }

\newcommand{\la}{\lambda }
\newcommand{\mi}{\mu }
\newcommand{\ks}{\xi}

\newcommand{\ti}{\tau }

\newcommand{\C}{\mathbb{C}}

\newcommand{\R}{\mathbb{R}}

\newcommand{\N}{\mathbb{N}}

\newcommand{\K}{\mathbb{K}}

\newcommand{\ssum}{\sum\limits}

\newcommand{\tl}{\widetilde{L}}

\newcommand{\ld}{\ldots}

\newcommand{\qb}{$\quad\blacksquare$}

\begin{document}
\pagestyle{myheadings}
\markboth{Universal Taylor Series on products of planar domains}{K. Kioulafa, G. Kotsovolis, V. Nestoridis}
\title{\bf Universal Taylor Series on products of planar domains}
%
%
\author{K. Kioulafa, G. Kotsovolis, V. Nestoridis}
\date{}
\maketitle
\begin{abstract}
Using a recent Mergelyan type theorem for products of planar compact sets we establish generic existence of Universal Taylor Series on products of planar simply connected domains $\OO_i$, $i=1,\ld,d$. The universal approximation is realized by partial sums of the Taylor development of the universal function on products of planar compact sets $K_i$, $i=1,\ld,d$ such that $\C-K_i$ is connected and for at least one $i_0$ the set $K_{i_0}$ is disjoint from $\OO_{i_0}$.
\end{abstract}
{\em AMS classification numbers}: 30K05, 32A05, 32A17, 30A30\smallskip
\section{Introduction}\label{sec1}
\noindent

The proof of the existence of universal Taylor series in one complex variable (\cite{CP}, \cite{Luh}, \cite{Nes1}), \cite{Me-Ne} and \cite{Ka}, \cite{G-E}, \cite{Nes2}, \cite{K-K-N} is based on Mergelyan's theorem \cite{Ru}. In several complex variables approximation theory is less developed and there is no satisfactory Mergelyan's theorem. In consequence there are only a few theorems for the existence of some kind of universal Taylor series in several variables (\cite{Cl}, \cite{Da-Ne}). Recently, during a Research in pairs program at the Institute of Oberwolfach \cite{F.G.M.N} a new Mergelyan's type theorem was obtained for products of compact planar sets. In particular, if $K_1,\ld,K_d$ are compact subsets of the plane $\C$ with connected complements $\C-K_i$, $i=1,\ld,d$, then there is a Mergelyan's type theorem for the product $K=K_1\times\cdots\times K_d$. Every function $f:K\ra\C$ continuous on $K$ and such that $f\circ\Fi$ is holomorphic on the disk $D\subset\C$ for every injective holomorphic mapping $\Fi:D\ra K$, is uniformly approximated on $K$ by polynomials of $d$ complex variables. Using the above fact, we prove in the present paper generic universality of Taylor series on products of planar simply connected domains $\OO_i$, $i=1,\ld,d$, where the universal approximation is valid on products of planar compact sets $K_i$, $i=1,\ld,d$ where $\C-K_i$ are connected and for at least one $i_0$ the set $K_{i_0}$ is disjoint from $\OO_{i_0}$. The approximation is realized by partial sums of the Taylor development of the generic universal function on $\OO=\prod\limits^d_{i=1}\OO_i$. There are a lot of ways to enumerate all the monomials of a Taylor expansion in severable variables. We prove that for any such enumeration there exist Universal Taylor Series and their set is $G_\de$ and dense in the space $H(\OO)$ of holomorphic functions on $\OO$ endowed with its natural Fr\'{e}chet topology. Choosing appropriately the above enumeration we can obtain the universal approximation using spherical partial sums $(n^2_1+\cdots+n^2_d\le N^2)$ or partial sums in a rectangular shape $(\max(n_1,\ld,n_d)\le N)$ etc. We show that there is no holomorphic function on $\OO$, which is universal with respect to all enumerations.

The method of our proof is by using the abstract theory of universal series (\cite{Ne-Pa}, \cite{B-Gr-Erd-N-P}) which is based on Baire's Category theorem. For the use of Baire's theorem in analysis we refer to \cite{Ka}, \cite{G-E}.
\section{Abstract Theory of universal series}\label{sec2}
\noindent

We now introduce some context from \cite{B-Gr-Erd-N-P}, which will be of significant help towards our results. Let us assume that $(X_\ti)$, $\ti\ge1$ is a sequence of metrizable topological vector spaces over the field $\K$ $(\K=\C$ or $\K=\R)$, equipped with translation - invariant metrics $\rho_\ti$. Denote by $E$ a complete metrizable topological vector space, with topology induced by a translation - invariant metric $d$. Let $L_m$, $m\ge1$ be an increasing sequence of non-empty compact subsets of the metric space $L=\bigcup\limits_{m\ge1}L_m$ and $\ks_0$ a distinguished element of $L_1$. Suppose that for any $n\in\N$, we have continuous maps
\[
e_n:L\ra E, \ \ x_{\ti,n}:L\ra X_\ti \ \ \text{and} \ \ \f_n:L\times E\ra\K.
\]
Let $G$ denote the set of sequences in $\K^N$ with finitely many non-zero elements and let us make the following assumptions:
\begin{enumerate}
\item[i)] The set $\{g_a:a\in G\}$ is dense in $E$ and $g_a=\ssum^\infty_{j=0}a_je_j(\ks_0)$, where $a=\big(a_j\big)^\infty_{j=0}$.
\item[ii)] For every $a\in G$ and every $\ks\in L$, the sets $\{n:\f_n(\ks,g_a)\neq0\}$ are finite and uniformly bounded with respect to $\ks\in L_m$ for any $m\ge1$.
\item[iii)] For every $a\in G$ and $\ks\in L$:
\[
\sum^\infty_{j=0}\f_j(\ks,g_a)e_j(\ks)=\sum^\infty_{j=0}a_je_j(\ks_0).
\]
\item[iv)] For every $a\in G$, $\ks\in L$ and $\ti\ge 1$:
\[
\sum^\infty_{j=0}\f_j(\ks,g_a)x_{\ti,j}(\ks)=\sum^\infty_{j=0}a_jx_{\ti,j}(\ks_0).
\]
\end{enumerate}
\begin{Def}\label{def2.1}
Let $\mi$ be an infinite subset of $\{0,1,2,\ld\}$. Under the above assumptions, an element $f\in E$ belongs to the class $U^\mi_{E,L}$, if, for every $\ti\ge1$ and every $x\in X_\ti$, there exists a sequence $\la=(\la_n)\subset\mi$, such that, for every $m\ge1$ we have
\begin{enumerate}
\item[i)] $\sup\limits_{\ks\in L_m}\rho_\ti\Big(\ssum^{\la_n}_{j=0}\f_j(\ks,f)x_{\ti,j}(\ks),x\Big)\ra0$ as $n\ra+\infty$, for every $\ti\ge1$.
\item[ii)] $\sup\limits_{\ks\in L_m}d\Big(\ssum^{\la_n}_{j=0}\f_j(\ks,f)e_j(\ks),f\Big)\ra0$ as $n\ra+\infty$.
\end{enumerate}
In the case $\mi=\{0,1,2,\ld\}$ we write $U_{E,L}$ instead of $U^\mi_{E,L}$.
\end{Def}
\begin{rem}\label{rem2.2}
In Definition \ref{def2.1} it is equivalent to require that the sequence $\la=(\la_n)\subset\mi$ is strictly increasing \cite{Ve}.
\end{rem}
\begin{thm}\label{thm2.3} \cite{B-Gr-Erd-N-P}.\\
Under the above assumptions, the following are equivalent.

1) For every $\ti\ge1$, every $x\in X_\ti$ and every $\e>0$, there exists $n\in N$ and $a_0,a_1,\ld,a_n\in\K$ so that
\[
\rho_\ti\bigg(\sum^n_{j=0}a_jx_{\ti,j}(\ks_0),x\bigg)<\e \ \ \text{and} \ \ d\bigg(\sum^n_{j=0}a_je_j(\ks_0),0\bigg)<\e.
\]

2)For any increasing sequence $\mi$ of positive integers the set $U^\mi_{E,L}$ is a dense - $G_\de$ subset of $E$ and $U^\mi_{E,L}\cup\{0\}$ contains a dense subspace of $E$.
\end{thm}
\section{The function Algebra $\bbb{A_D(K)}$ and two approximation lemmas}\label{sec3}
\noindent

In this section we present some new results from \cite{F.G.M.N} which will be useful later on.
\begin{Def}\label{def3.1} \cite{F.G.M.N}\\
Recalling that a mapping from a planar domain to $C^d$ is holomorphic, if each coordinate is a complex valued holomorphic function, we define the algebra $A_D(K)$ as the set of all functions $f:K\ra C$, where $K$ is compact in $C^d$, that are continuous on $K$ and such that, for every disc $D\subset\C$ and every holomorphic mapping $\f:D\ra K\subset C^d$, the composition $f\circ\f:D\ra C$ is holomorphic.
\end{Def}
\begin{lem}\label{lem3.2} \cite{F.G.M.N}\\
Let $K=\prod\limits^d_{i=1}K_i$, where $K_i\subset C$ is compact, with $C-K_i$ connected. If $f\in A_D(K)$ and $\e>0$, then there exists a polynomial $P$ such that $\sup\limits_{z\in K}|f(z)-P(z)|<\e$.
\end{lem}

Lemma \ref{lem3.2} is part of a result in \cite{F.G.M.N}, Remark 4.9 (3). Furthermore, from the same paper, we also have the following lemma (See Prop. 2.8 of \cite{F.G.M.N}).

\begin{lem}\label{lem3.3} \cite{F.G.M.N}\\
Let $K$ be as in Lemma \ref{lem3.2}. If $f$ is holomorphic on a neighborhood of $K$, $I$ is a finite subset of $N^d$ and $\e>0$, there exists a polynomial $Q$ such that:
\[
\sup_{z\in K}\bigg|\frac{\partial^{a_1+\cdots+a_d}}{\partial ^{a_1}_{z_1}\partial^{a_2}_{z_2}\cdots\partial^{a_d}_{z_d}}(f-Q)(z)\bigg|<\e
\]
for every $a=(a_1,a_2,\ld,a_d)\in I$.
\end{lem}
\section{$\bbb{H(\OO)}$ universal Taylor series on products of planar domains}\label{sec4}
\noindent

Let $\OO\subseteq C^d$ be an open set, $\zi=(\zi_1,\zi_2,\ld,\zi_d)$ a point of $\OO$ and $N_j$, $j=0,1,2,\ld$ be an enumeration of $N^d$, where $N=\{0,1,2,\ld\}$, and denote the monomials by $(z-\zi)^a=(z_1-\zi_1)^{a_1}\cdots(z_d-\zi_d)^{a_d}$, where $a=(a_1,a_2,\ld,a_d)$. For $f$ a holomorphic function on $\OO$, we denote by $a(f,\zi)$ the coefficient of $(z-\zi)^a$ in the Taylor expansion of $f$ with center $\zi$; that is,
\[
a(f,\zi)=\frac{1}{a_1!a_2!\cdots a_d!}\frac{\partial^a}{\partial^{a_1}_{z_1}\cdots\partial^{a_d}_{z_d}}f(\zi).
\]
Let
\[
S_N(f,\zi)(z)=\sum^N_{j=0}a_{N_j}(f,\zi)(z-\zi)^{N_j}.
\]
\begin{thm}\label{thm4.1}
Let $\OO_i\subseteq C$ simply connected domains, $i=1,2,\ld,d$ and $\OO=\prod\limits^d_{i=1}\OO_i$. Fix $\zi_0\in\OO$ and the enumeration $N_j$ as above. Let $\mi$ be an infinite subset of $N$. There exists a holomorphic function $f\in H(\OO)$ such that for every compact sets $K_i\subseteq C-\OO_i$ with $C-K_i$ connected and every $h\in A_D\Big(\prod\limits^d_{i=1}K_i\Big)$, there exists a strictly increasing sequence $\la_n\in\mi$, $n=1,2,\ld$ such that
\[
\begin{array}{l}
  S_{\la_n}(f,\zi^0)(z)\ra h(z) \ \ \text{uniformly on} \ \ \prod\limits^d_{i=1}K_i \ \ \text{and} \\
  S_{\la_n}(f,\zi^0)(z)\ra f(z) \ \ \text{uniformly on compacta of $\OO$ as $n\ra+\infty$}.
\end{array}
\]
Furthermore, the set of such $f\in H(\OO)$ is a dense-$G\de$ subset of the space $H(\OO)$ of holomorphic functions on $\OO$ endowed with the topology of uniform convergence on compacta and contains a dense vector subspace except 0.
\end{thm}
\begin{Proof}
According to Remark \ref{rem2.2} it is not necessary to prove $\la_n<\la_{n+1}$; it suffices that $\la_n\in\mi$ for all $n$.

It is known \cite{Nes2}, \cite{Me-Ne}, that when $\OO_i$ is a (simply connected) domain in $C$, there can be found compact sets $K_{i,1},K_{i,2},K_{i,3},\ld$ in $C-\OO_i$, with $C-K_{i,j}$ connected, such that for every $K$ compact set in $C-\OO_i$, with $C-K$ connected there exists $j:K\subseteq K_{i.j}$. Let $K_\ti$, $\ti=1,2,\ld$ be an enumeration of all $\prod\limits^d_{i=1}K_i,j_i$ where $j_i\in\{1,2,\ld\}$. We notice that, if we have the theorem for compact sets $K$ of the form $K_\ti$, we also have it for every $K=\prod\limits^d_{i=1}T_i$, where $T_i$ are compact sets in $C-\OO_i$, with $C-T_i$ connected. Indeed let us suppose that the function $f\in H(\OO)$ satisfies the theorem for compact sets $K$ of the form $K_\ti$ and $h\in A_D\Big(\prod\limits^d_{i=1}T_i\Big)$. Let $\e>0$. By Lemma \ref{lem3.2}, there exists a polynomial $P$ such that $\dis\sup_{z\in\prod\limits^d_{i=1}T_i}|P(z)-h(z)|<\dfrac{\e}{2}$. By definition of the sets $K_\ti$, we can also find $\ti\in N$, such that $\prod\limits^d_{i=1}T_i\subseteq K_\ti$. If $d$ is the standard metric in $H(\OO)$, inducing the topology of uniform convergence on compacta of $\OO$, there exists $\la\in \mi$ such that
\[
d(S_\la(f,\zi^0),f)<\e \ \ \text{and} \ \ \sup_{z\in K_\ti}|S_\la(f,\zi^0)(z)-P(z)|<\frac{\e}{2}.
\]
The result follows by the triangular inequality.\\
We proceed now to prove the result for the compact sets $K_\ti$, $\ti=1,2,\ld$. Following the notation of Section \ref{sec2}, we set $X_\ti=A_D(K_\ti)$, $E=H(\OO)$, $L_m=\{\zi_0\}$ for every $m\in N$ and $\ks_0=\zi_0$. Take $e_n(\ks)=x_{\ti,n}(\ks)=(z\ra(z-\ks)^{N_n})$ and $\f_n(\ks,f)=a_{N_n}(f,\ks)$. One can easily check that conditions i), ii), iii) and iv) of \S\,2 are verified. Thus, according to Theorem \ref{thm2.3}, it suffices to show that, for every $\ti\ge1$, $g\in A_D(K_\ti)$, every $\e>0$, and $\widetilde{L}\subseteq\OO$ compact, there exists a polynomial $P$ such that
\[
\begin{array}{l}
  \sup\limits_{z\in K_\ti}|P(z)-g(z)|<\e, \ \ \text{and}\\
  \sup_{z\in\widetilde{L}}|P(z)|<\e.
\end{array}
\]
Since $\OO$ is a product of simply connected domains $\OO_i\subset C$, without loss of generality, we can assume that $\tl=\prod\limits^d_{i=1}\tl_i$, where $\tl_i\subset\OO_i$ is compact and $C-\tl_i$ is connected \cite{Ru}.

Since $g\in A_D(K_\ti)$, according to Lemma \ref{lem3.2}, there exists a polynomial $\widetilde{g}$, such that $\sup\limits_{z\in K_\ti}|g(z)-\widetilde{g}(z)|<\dfrac{\e}{2}$. Denote also $K_\ti=\prod\limits^d_{i=1}K_{i,j_i}$. We define the function $h:\prod\limits^d_{i=1}(\tl_i\cup K_{i,j_i})\ra C$ as follow:
\[
\begin{array}{l}
  h(z)=\widetilde{g}(z) \ \ \text{for} \ \ z\in\prod\limits^d_{i=1}K_{i,j_i} \\
  h(z)=0 \ \ \text{for} \ \ z\in\Big[\prod\limits^d_{i=1}(\tl_i\cup K_{i,j_i})\Big]\backslash\prod\limits^d_{i=1}K_{i,j_i}
\end{array}
\]
We notice that $\widetilde{g}$, being a polynomial, is everywhere defined. Therefore, other functions $h$ could be defined, as well.\\
Since the function $h$ is holomorphic on a neighborhood of $\prod\limits^d_{i=1}(K_{i,j_i}\cup\tl_i)$, it follows easily that $h\in A_D\Big(\prod\limits^d_{i=1}(\tl_i\cup K_{i,j_i})\Big)$.

The result follows from Lemma \ref{lem3.2} applied for $\dfrac{\e}{2}$ and the triangular inequality. \qb
\end{Proof}
\begin{thm}\label{thm4.2}
Under the assumptions and notation of Theorem \ref{thm4.1}, there exists a holomorphic function $f\in H(\OO)$, such that for every compact sets $K_i\subseteq C-\OO_i$ with $C-K_i$ connected and every $h\in A_D\Big(\prod\limits^d_{i=1}K_i\Big)$, there exists a strictly increasing sequence $\la_n\in\mi$ such that for every compact set $\tl\subset\OO$ we have
\[
\begin{array}{l}
  \sup\limits_{z\in\prod\limits^d_{i=1}K_i,\zi\in\tl}|S_{\la_n}(f,\zi)(z)-h(z)|\ra0 \ \ \text{and} \\
  \sup\limits_{\zi\in \widetilde{L},z\in\tl}|S_{\la_n}(f,\zi)(z)-f(z)|\ra0 \ \ \text{as} \ \ n\ra+\infty.
\end{array}
\]
Furthermore, the set of such functions $f\in H(\OO)$ is a dense - $G\de$ subset and contains a dense vector subspace except 0.
\end{thm}
\begin{Proof}
The proof is similar to the proof of Theorem \ref{thm4.1}. The only difference is that when we apply the abstract theory of \S\,2 we do not set $L_m=L=\{\zi_0\}$ but we set $L=\OO$, $L_m$, $m=1,2,\ld$ an exhausting family of compact subsets of $\OO$, where each $L_m$ is a product of planar compact sets with connected complements and $\zi_0$ an arbitrary element of $L_1$. \qb
\end{Proof}

For the next two theorems we will need a new definition.
\begin{Def}\label{def4.3}
Let $K$ be a compact set of $C^d$. We note by $O(K)$, the set of functions $g$ holomorphic on an open set $V_g$ such that $K\subset V_g\subset C^d$. The topology of the space $O(K)$ is the Fr\'{e}chet topology defined by the denumerable family of seminorms $\sup\limits_{\zi\in K}|D_a(\zi)|$, where $D_a$ varies in the denumerable set of all differential operators of  mixed partial derivatives in $\zi=(\zi_1,\zi_2,\ld,\zi_d)$. We denote by $\rho$ the corresponding metric.
\end{Def}
\begin{thm}\label{thm4.4}
Under the assumptions and notation of Theorem \ref{thm4.1}, there exists a holomorphic function $f\in H(\OO)$, such that, for every compact sets $K_i\subseteq C-\OO_i$, with $C-K_i$ connected and every function $h\in O\Big(\prod\limits^d_{i=1}K_i\Big)$, there exists a strictly increasing sequence $\la_n\in\mi$, $n=1,2,\ld$, such that for every differential operator $D$ of mixed partial derivatives with respect to $z=(z_1,z_2,\ld,z_d)$ we have
\[\begin{array}{l}
    D(S_{\la_n}(f,\zi^0))(z)\ra D(h)(z) \ \ \text{uniformly on} \ \ \prod\limits^d_{i=1}K_i \ \ \text{and} \\
    S_{\la_n}(f,\zi^0)(z)\ra f(z) \ \ \text{uniformly on compacta of $\OO$ as $n\ra+\infty$}.
  \end{array}
\]
Furthermore, the set of such $f\in H(\OO)$ is a dense - $G\de$ subset of the space $H(\OO)$ and contains a dense vector subspace except 0.
\end{thm}
\begin{Proof}
We keep the notation of the compact sets $K_\ti$, as in the proof of Theorem \ref{thm4.1}. We notice, that if we have the theorem for compact sets $K$ of the form $K_\ti$, we also have it for every $K=\prod\limits^d_{i=1}T_i$, where $T_i$ is compact in $C-\OO_i$, with $C-T_i$ connected. Indeed let us assume that the function $f\in H(\OO)$ satisfies the theorem for compact sets $K$ of the form $K_\ti$ and let $h\in O\Big(\prod\limits^d_{i=1}T_i\Big)$. Let $\e>0$. By Lemma \ref{lem3.3} there exists a polynomial $P$ such that $\rho(P,h)<\dfrac{\e}{2}$, where $\rho$ is the metric in $O\Big(\prod\limits^d_{i=1}T_i\Big)$, (Definition \ref{def4.3}). By the definition of the sets $K_\ti$, we can also find $\ti\in N$, such that $\prod\limits^d_{j=1}T_j\subseteq K_\ti$. If $d$ is the standard metric in $H(\OO)$, inducing the topology of uniform convergence on compacta of $\OO$, there exists $\la\in \mi$ such that
\[
\begin{array}{l}
  d(S_\la(f,\zi^0),f)<\e \ \ \text{and} \\ [1.5ex]
  \rho(S_\la(f,\zi^0),P)<\dfrac{\e}{2}.
\end{array}
\]
The result follows by the triangular inequality. \\
We proceed now to prove the result for the compact sets $K_\ti$, $\ti=1,2,\ld$\;.
Following the notation of Section \ref{sec2}, we set $X_\ti=O(K_\ti)$, $E=H(\OO)$, $L_m=\{\zi_0\}$ for every $m\in N$ and $\ks_0=\zi_0$. One can easily check that conditions i), ii), iii) and iv) of \S\,2 are verified. Thus, according to Theorem \ref{thm2.3} it suffices to show for every $\ti\ge1$, $g\in O(K_\ti)$, every $\e>0$, and $\tl\subseteq\OO$ compact, there exists a polynomial $P$ such that
\[
\rho(P,g)<\e \ \ \text{and} \ \ \sup_{z\in\tl}|P(z)|<\e, \ \ \text{where $\rho$ is the metric of $O(K_\ti)$}.
\]
Since $\OO$ is a product of simply connected domains $\OO_i\subset C$, without loss of generality, we can assume that $\tl=\prod\limits^d_{i=1}\tl_i$, where $\tl_i\subset\OO_i$ is compact and $C-\tl_i$ is connected \cite{Ru}.

Since $g\in O(K_\ti)$, according to Lemma \ref{lem3.3}, there exists a polynomial $\widetilde{g}$, such that $\rho(g,\widetilde{g})<\dfrac{\e}{2}$. Denote also $K_\ti=\prod\limits^d_{i=1}K_{i,j_i}$. We define the function
\[
\begin{array}{l}
h:\prod\limits^d_{i=1}(\tl_i\cup K_{i,j_i})\ra C \ \ \text{as follows}   \\
  h(z)=\widetilde{g}(z) \ \ \text{for} \ \ z\in\prod\limits^d_{i=1}K_{i,j_i} \\
  h(z)=0 \ \ \text{for} \ \ z\in\Big[\prod\limits^d_{i=1}(\tl_i\cup K_{i,j_i})\Big]\backslash\prod\limits^d_{i=1}K_{i,j_i}.
\end{array}
\]

We notice that $\widetilde{g}$ is everywhere defined, because it is a polynomial. Therefore, other functions $h$ could also be defined, as well.

Since the function $h$ is holomorphic on a neighborhood of $\prod\limits^d_{i=1}(\tl_i\cup K_{i,j_i})$, it follows by Lemma \ref{lem3.3}, that there exists a polynomial $P$ such that $\rho'(h,P)<\dfrac{\e}{2}$ where $\rho'$ is the metric of $O\Big(\prod\limits^d_{i=1}(K_{i,j_i}\cup\tl_i)\Big)$. Thus, we also have
\[
\rho(h,P)<\frac{\e}{2} \ \ \text{and} \ \ \sup_{z\in\tl}|P(z)-h(z)|<\frac{\e}{2}.
\]
The result follows by the triangular inequality. \qb
\end{Proof}
\begin{thm}\label{thm4.5}
Under the assumptions and notation of Theorem \ref{thm4.2}, there exists a holomorphic function $f\in H(\OO)$, such that for every compact sets $K_i\subseteq C-\OO_i$, with $C-K_i$ connected and every function $h\in O\Big(\prod\limits^d_{i=1}K_i\Big)$, there exists a strictly increasing sequence $\la_n\in\mi$, $n=1,2,\ld$, such that for every differential operator $D$ of mixed partial derivatives with respect to $z=(z_1,z_2,\ld,z_d)$ and for every compact set $\tl\subseteq\OO$ we have that
\[
\begin{array}{l}
  \sup\limits_{\zi\in\tl,z\in\prod\limits^d_{i=1}K_i}|D(S\la_n(f,\zi))
  (z)-D(h)(z)|\ra0 \ \ \text{and}  \\
  \sup\limits_{\zi\in\tl,z\in\tl}|S\la_n(f,\zi)(z)-f(z)|\ra0 \ \ \text{and} \ \ \text{as} \ \ n\ra+\infty.
\end{array}
\]
Furthermore, the set of such functions $f\in H(\OO)$ is a dense - $G\de$ subset of the space $H(\OO)$ and contains a dense vector subspace except 0.
\end{thm}
\begin{Proof}
The proof is similar to the proof of Theorem \ref{thm4.4}. The only difference is that we do not set $L_m=L=\{\zi_0\}$ but we set $L=\OO$, $L_m$, $m=1,2,\ld$ an exhausting family of compact subsets of $\OO$, where each $L_m$ is a product of planar compact sets with connected complements and $\zi_0$ an arbitrary element of $L_1$. \qb
\end{Proof}
\begin{rem}\label{rem4.6}
We note that the previous theorems have been shown for products of compact sets $K_i$, where $C-K_i$ are connected and $K_i$ is contained in $C-\OO_i$. However, without any particular effort, the previous theorems are also valid for products of compact sets $K_i$, where $C-K_i$ are connected and there exists at least one $i_0$, such that $K_{i_0}$ is contained in $C-\OO_{i_0}$. This is equivalent to say that $C-K_i$ is connected and the product of $K_i$ is disjoint from the product of $\OO_i$. The proof in this case is the same , except from the part where the function $h$ is defined. For every $i\neq i_0$ there exists a disc $B_i$ in $\C$, such that $L_i\cup K_i$ is contained in $B_i$.The function $h$ in this case should be defined as a function $h: \prod\limits S_i\ra C$, where $S_i=B_i$ for $i\neq i_0$ and $S_{i_0}=L_{i_0}\cup K_{i_0}$. We set $h(z)=\widetilde{g}(z)$ for $z$ in the product of $B_i$ for $i\neq i_0$ and $K_{i_0}$, and $h(z)=0$ for $z$ in the product of $B_i$ for $i\neq i_0$ and $L_{i_0}$. We notice that $\widetilde{g}$, being a polynomial, is everywhere defined.
\end{rem}
\section{Smooth universal Taylor series}\label{sec5}
\noindent

For $\OO\subseteq C^d$ open, we denote $A^\infty(\OO)$, the set of all holomorphic functions $f\in H(\OO)$, such that for every differential operator $D$ of mixed partial derivatives in $z=(z_1,z_2,\ld,z_d)$, the function $Df$ has a continuous extension on $\overline{\OO}$. The topology of the space $A^\infty(\OO)$ is the Fr\'{e}chet topology defined by the denumerable family of seminorms $\sup\limits_{\zi\in\overline{\OO},|\zi|\le n}|Df(\zi)|$, where $n\in\{1,2,\ld\}$ and $D$ varies in the denumerable set of all differential operators of mixed partial derivatives in $\zi=(\zi_1,\zi_2,\ld,\zi_d)$. We notice that for $\zi\in\partial\OO$, $Df(\zi)$ is defined by continuity of $Df$ in $\overline{\OO}$.

Let $P$ be the set of all polynoomials in $d$ variables and $X^\infty(\OO)$ be the closure of $P$ in $A^\infty(\OO)$.
\begin{thm}\label{thm5.1}
Under the assumptions of Theorem \ref{thm4.1}, we also assume that $C-\overline{\OO}_i$ is connected for each $i=1,\ld,d$. Then there exists a holomorphic function $f\in X^\infty(\OO)$, such that, for every compact sets $K_i\subset C-\overline{\OO}_i$, with $C-K_i$ connected and every $h\in A_D\Big(\prod\limits^d_{i=1}K_i\Big)$, there exists a strictly increasing sequence $\la_n\in\mi$, $n=1,2,\ld$ such that
\[
\begin{array}{l}
  S\la_n(f,\zi^0)(z)\ra h(z) \ \ \text{uniformly on} \ \ \prod\limits^d_{i=1}K_i \ \ \text{and} \\
  S\la_n(f,\zi^0)(z)\ra f(z) \ \ \text{in the topology of} \ \ A^\infty(\OO) \ \ \text{as} \ \ n\ra+\infty.
\end{array}
\]
Furthermore, the set of such functions $f\in X^\infty(\OO)$ is a dense - $G\de$ subset and contains a dense vector subspace except 0.
\end{thm}
\begin{Proof}
It is known \cite{CP}, \cite{Luh2} and \cite{K-K-N} that there exists a sequence of compact sets $K_{ij}$, $j=1,2,\ld$ in $C-\overline{\OO}_i$, such that $C-K_{i,j}$ is connected and for every $K$ compact in $C-\overline{\OO}_i$, with $C-K$ connected, there exists $j\ge1$ such that $K\subseteq K_{i,j}$. Let $K_\ti$ be an enumeration of all $\prod\limits^d_{i=1}K_{i,j_i}$. We notice, that if we have the theorem for the compact sets  $K$ of the form $K_\ti$, we also have it for every $K=\prod\limits^d_{j=1}T_j$, where $T_i$ is compact in $C-\overline{\OO}_i$, with $C-T_i$ connected. Indeed, let $f\in X^\infty(\OO)$ satisfy the theorem for the compact sets $K$ of the form $K_\ti$ and $h\in A_D\Big(\prod\limits^d_{j=1}T_j\Big)$. Let $\e>0$. By Lemma \ref{lem3.2}, there exists a polynomial $P$ such that $\sup\limits_{z\in\prod\limits^d_{j=1}T_j}|P(z)-h(z)|<\dfrac{\e}{2}$. By definition of the sets $K_\ti$, we can also find $\ti\in N$, such that $\prod\limits^d_{j=1}T_j\subseteq K_\ti$. If $d$ is the metric of $A^\infty(\OO)$, there exists $\la\in \mi$ such that
\[
\begin{array}{l}
  d(S_\la(f,\zi^0),f)<\e \ \ \text{and} \\  [1.5ex]
  \sup\limits_{z\in K_\ti}|S_\la(f,\zi^0)(z)-P(z)|<\dfrac{\e}{2}.
\end{array}
\]
The result follows by the triangular inequality. We proceed now to prove the result for the compact sets $K_\ti$.

Following the notation of Section \ref{sec2}, we set $X_\ti=A_D(K_\ti)$, $E=X^\infty(\OO)$, $L_m=\{\zi_0\}$ for every $m\in N$ and $\ks_0=\zi_0$. One can easily check that conditions i), ii), iii) and iv) of \S\,2 are verified. Thus, according to Theorem \ref{thm2.3}, it suffices to show that for every $\ti\ge1$, $g\in A_D(K_\ti)$, every $\e>0$, every finite set $F$ of mixed partial derivatives  in $z=(z_1,z_2,\ld,z_d)$ and every finite sequence $S=\{n_\el\}\subset\N$, there exists a polynomial $P$ such that
\[
\begin{array}{l}
  \sup\limits_{z\in K_\ti}|P(z)-g(z)|<\e \ \ \text{and} \\ [1ex]
  \sup\limits_{\zi\in\overline{\OO},|\zi|\le n_\el}|D_{\el'}(P)(\zi)|<\e \ \ \text{for} \ \ D_{\el'}\in F, \ \ n_\el\in S.
\end{array}
\]
Let $N$ be the maximum of the finite sequence $n_\el$. It obviously suffices to show that if $\tl=\overline{\OO}\cap\{z\in C^d:|z|\le N\}$, then
\[
\sup_{z\in K_\ti}|P(z)-g(z)|<\e \ \ \text{and} \ \ \sup_{z\in\tl}
|D_{\el'}(P)(z)|<\e \ \ \text{for} \ \ D_{\el'}\in F.
\]
Since $g\in A_D(K_\ti)$, there exists, by Lemma 3.2, a polynomial $\widetilde{g}$, such that $\sup\limits_{z\in K_\ti}|\widetilde{g}(z)-g(z)|<\dfrac{\e}{2}$. Moreover, since $\overline{\OO}=\prod\limits^d_{j=1}\overline{\OO}_j$, without loss of generality, we can assume that $\tl=\prod\limits^d_{j=1}\tl_j$, where $\tl_j\subset\overline{\OO}_j$ is compact and $C-\tl_j$ is connected. Denote also $K_\ti=\prod\limits^d_{i=1}K_{i,j_i}$. We define the function $h:\prod\limits^d_{i=1}(\tl_i\cup K_{i,j_i})\ra C$ as follows:
\[
\begin{array}{l}
  h(z)=\widetilde{g}(z) \ \ \text{for} \ \ z\in\prod\limits^d_{i=1}K_{i,j_i} \\ [1.5ex]
  h(z)=0 \ \ \text{for} \ \ z\in\Big[\prod\limits^d_{i=1}(\tl_i\cup K_{i,j_i})]\backslash\prod\limits^d_{i=1}K_{i,j_i}.
\end{array}
\]
We notice that $\widetilde{g}$, being a polynomial, is defined everywhere. Thus, other definitions of the function $h$ would also be possible.

Since the function $h$ is holomorphic on a neighborhood of $\prod\limits^d_{i=1}(\tl_i\cup K_{i,j_i})$, it follows by Lemma \ref{lem3.3}, that there exists a polynomial $P$ such that
\[
\begin{array}{l}
  \sup\limits_{z\in\prod\limits^d_{i=1}(\tl_i\cup  K_{i,j_i})}|P(z)-h(z)|<\dfrac{\e}{2} \ \ \text{and} \\
  \sup\limits_{z\in\prod\limits^d_{i=1}(\tl_i\cup K_{i,j_i})}|D_{\el'}(P)(z)-D_{\el'}(h)(z)|<\dfrac{\e}{2} \ \ \text{for} \ \ D_{\el'}\in F.
\end{array}
\]
The result follows by the triangular inequality. \qb
\end{Proof}

We notice that for a function $f\in X^\infty(\OO)$ it make sense to write $S_\la(f,\zi)(z)$ even if $\zi\in\partial\OO$ because of the continuous extension of the partial derivatives of $f$. 
\begin{thm}\label{thm5.2}
Under the assumptions of Theorem \ref{thm5.1}, there exists a holomorphic function $f\in X^\infty(\OO)$, such that, for every compact sets $K_i\subseteq C-\overline{\OO}_i$, with $C-K_i$ connected and every $h\in A_D\Big(\prod\limits^d_{i=1}K_i\Big)$, there exists a strictly increasing sequence $\la_n\in\mi$, $n=1,2,\ld$, such that, for every compact set $\tl\subset\overline{\OO}$ we have
\[
\begin{array}{l}
  \sup\limits_{z\in\prod\limits^d_{i=1}K_i,\zi\in\tl}|S{\la_n}(f,\zi)(z)-h(z)|\ra0 \ \ \text{and} \\
  \sup\limits_{\zi\in\tl}d(S_{\la_n}(f,\zi),f)\ra0 \ \ \text{as} \ \ n\ra+\infty, \ \ \text{where $d$ is the metric of $A^\infty(\OO)$}.
\end{array}
\]
Furthermore, the set of such functions $f\in X^\infty(\OO)$ is a dense - $G\de$ subset and contains a dense vector subspace except 0.
\end{thm}
\begin{Proof}
The proof is similar to the proof of Theorem \ref{thm5.1}. The only difference is that we do not set $L_m=L=\{\zi_0\}$ but we set $L=\overline{\OO}$, $L_m=\prod\limits^d_{i=1}L_{i,m}$ where $L_{i,m}=\overline{\OO}_i\cap\{z\in\C:|z|\le m\}$ and let $\zi_0$ be an arbitrary point of $L_1$. \qb
\end{Proof}
\begin{rem}\label{rem5.3}
A sufficient condition for the simply connected domains $\OO_i$ satisfying that $C-\overline{\OO}_i$ is connected, so that $X^\infty(\OO)=A^\infty(\OO)$ is that there exists $M_i<+\infty$ so that, for every $A,B\in\OO_i$, there exists a curve $\ga$ in $\OO_i$ joining $A$ and $B$ with length $(\ga)\le M_i$ \cite{Ne-Z}; see also \cite{S-S-N} and \cite{G-N100}.
\end{rem}
\begin{rem}\label{rem5.4}
We note that the previous theorems have been shown for products of compact sets $K_i$, where $C-K_i$ are connected and $K_i$ is contained in $C-\overline{\OO}_i$. However, without any particular effort , the previous theorems are also valid for products of compact sets $K_i$, where $C- K_i$ are connected and there exists at least one $i_0$, such that $K_{i_0}$ is contained in $\C-\overline{\OO}_{i_0}$. The alteration of the proof is similar to that of Remark \ref{rem4.6}.
\end{rem}
\section{Universal Taylor series with respect to several enumerations}\label{sec6}
\noindent

It is an obvious fact that due to Baire's theorem , for every denumerable family of enumerations of the set of mononyms ,the set of unviersal series with respect to all enumerations of this family is a $G\de$ – dense set. From this arises the question whether or not there exists a universal series with respect to all possible enumerations of the mononyms.

The following theorem provides a negative answer to this question.
\begin{thm}\label{thm6.1}
With the notation of Theorem \ref{thm4.1}, for every holomorphic function $f\in H(\OO)$ there exists an enumeration $N''_j$, such that the function f is not universal with respect to the enumeration $N''_j$.
\end{thm}
\begin{Proof}
Let $z_1$ be an element of the product $\OO$ of the sets $\OO_i$ and $z_2$ outside $\OO$. Define $A_j$, $j\ge0$ the sequence where $A_j=a_{N_j}(f,z_1)(z_2-z_1)^{N_j}$ for any arbitrary enumeration $N_j$. Let now, $R_j$ be the real part of $A_j$. The theorem is now proved if we take into consideration the following lemma. \qb
\end{Proof}
\begin{lem}\label{lem6.2}
Let $b_n$, $n\ge0$ a sequence of real numbers. Then, there exists a rearrangement of $\Sig b_n$ such that its partial sums are not dense in $R$.
\end{lem}

The proof of Lemma \ref{lem6.2} is elementary and is omitted. We can also prove that the partial sums of the rearrangement converges to some $\el$ in $[-\infty,+\infty]$.

We also notice that if we group together all the monomials of the same degree we have a similar result for the series of homogenous  terms.

\bigskip
\noindent
National and Kapodistrian University of Athens\\
Department of Mathematics\\
Panepistemiopolis, 157 84\\
Athens,
Greece \bigskip\\
e-mail addresses: \\
keiranna@math.uoa.gr \\
georgekotsovolis@yahoo.com\\
vnestor@math.uoa.gr \\

\end{document}